\documentclass[journal,twoside,web]{ieeecolor}
\usepackage{lcsys}
\usepackage{cite}
\usepackage{amsmath,amssymb,amsfonts}
\usepackage{algorithmic}
\usepackage{graphicx}
\usepackage{textcomp}
\def\BibTeX{{\rm B\kern-.05em{\sc i\kern-.025em b}\kern-.08em
    T\kern-.1667em\lower.7ex\hbox{E}\kern-.125emX}}
\usepackage[hidelinks]{hyperref}
\hypersetup{
    colorlinks=false,
}
\newtheorem{theorem}{Theorem}
\newtheorem{proposition}{Proposition}
\newtheorem{definition}{Definition}
\newtheorem{assumption}{Assumption}

\newcommand{\Nat}{{\mathbb{N}}}
\newcommand{\Real}{{\mathbb{R}}}
\newcommand{\N}{{\mathcal{N}}}
\newcommand{\T}{{\mathcal{T}}}
\newcommand{\X}{{\mathcal{X}}}
\newcommand{\Pop}{{\mathcal{P}}}
\newcommand{\D}{{\mathcal{D}}}
\newcommand{\A}{{\mathcal{A}}}
\newcommand{\G}{{\mathcal{G}}}
\newcommand{\Soc}{{\mathcal{S}}}
\newcommand{\sd}{d^*}
\newcommand{\epi}{\pi^*}
\newcommand{\esoc}{s^*}

\usepackage[acronym]{glossaries}
\newacronym{EV}{EV}{Electric Vehicle}
\newacronym{NE}{NE}{Nash Equilibrium}
\newacronym{NEs}{NE}{Nash Equilibria}
\newacronym{MFE}{MFE}{Mean Field Equilibrium}
\newacronym{SMFE}{SMFE}{Stationary Mean Field Equilibrium}
\newacronym{SNE}{SNE}{Stationary Nash Equilibrium}
\newacronym{SNEs}{SNE}{Stationary Nash Equilibria}
\newacronym{DPG}{DPG}{Dynamic Population Game}
\newacronym{DPGs}{DPGs}{Dynamic Population Games}
\newacronym{MDP}{MDP}{Markov Decision Process}
\newacronym{ODE}{ODE}{Ordinary Differential Equation}
\newacronym{MPE}{MPE}{Markov Perfect Equilibrium}

\title{Dynamic Population Games: A Tractable Intersection of Mean-Field Games and Population Games}

\usepackage[absolute,overlay]{textpos}
\usepackage[dvipsnames]{xcolor}

\begin{document}

\begin{textblock*}{\textwidth}(17mm,12mm) 
\noindent \bf \textcolor{NavyBlue}{Published in \emph{IEEE Control Systems Letters,} May 2024. \url{http://doi.org/10.1109/LCSYS.2024.3406947}}
\end{textblock*}

\begin{textblock*}{\textwidth}(17mm,69mm) 
\noindent \small \bf \textcolor{NavyBlue}{\textcopyright 2024 IEEE.  Personal use of this material is permitted.  Permission from IEEE must be obtained for all other uses, in any current or future media, including reprinting/republishing this material for advertising or promotional purposes, creating new collective works, for resale or redistribution to servers or lists, or reuse of any copyrighted component of this work in other works.}
\end{textblock*}

\author{Ezzat Elokda, Saverio Bolognani, \IEEEmembership{Member, IEEE}, Andrea Censi, Florian D\"orfler, \IEEEmembership{Senior Member, IEEE} and Emilio Frazzoli, \IEEEmembership{Fellow, IEEE}
\thanks{}
\thanks{Submitted on March 8, 2024. Research supported by NCCR Automation, a National Centre of
Competence in Research, funded by the Swiss National Science
Foundation (grant number 180545).}
\thanks{All authors are with ETH Zurich, 8092 Zurich, Switzerland.  (e-mails: \{elokdae,bsaverio,acensi,dorfler,efrazzoli\}@ethz.ch)}
}

\maketitle
\thispagestyle{empty}

\begin{abstract}

In many real-world large-scale decision problems, self-interested agents have individual dynamics and optimize their own long-term payoffs.
Important examples include 
the competitive access to shared resources (e.g., roads, energy, or bandwidth) but also non-engineering domains like epidemic propagation and control.
These problems are natural to model as \emph{mean-field games}.
Existing mathematical formulations of mean field games have had limited applicability in practice, since they require solving non-standard initial-terminal-value problems that are tractable only in limited special cases.
In this letter, we propose a novel formulation, along with computational tools, for a practically relevant class of \emph{\gls{DPGs}}, which correspond to discrete-time, finite-state-and-action, stationary mean-field games.
Our main contribution is a mathematical reduction of \emph{\gls{SNEs}} in \gls{DPGs} to standard \emph{\gls{NEs}} in static population games.
This reduction is leveraged to guarantee the existence of a \acrshort{SNE}, develop an \emph{evolutionary dynamics}-based \acrshort{SNE} computation algorithm, and derive simple conditions that guarantee stability and uniqueness of the \acrshort{SNE}.
We 
provide two examples of applications: fair resource allocation with heterogeneous agents and control of epidemic propagation.

Open source software for \acrshort{SNE} computation: \url{https://gitlab.ethz.ch/elokdae/dynamic-population-games}

\end{abstract}


\section{Introduction}
\label{sec:Intro}

\IEEEPARstart{M}{ean} field games, introduced under that name in~\cite{huang2006large,lasry2007mean}, and otherwise referred to as \emph{anonymous sequential games}~\cite{jovanovic1988anonymous,adlakha2015equilibria}, model competitive settings in large populations of selfish agents, who have individual state dynamics and optimize individual long-term average or discounted payoffs.
The microscopic state dynamics and payoffs of different agents are coupled through the \emph{mean field}, that is, the macroscopic distribution of states (and in some works, the distribution of states and actions~\cite{jovanovic1988anonymous,adlakha2015equilibria,guo2023general}).
Mean field games occur naturally in many important real-world problems.
For example, in an epidemic outbreak, agents have individual infection states, and the probability of getting infected depends on the proportion of infected agents, as well as whether the agents adhere to social distancing policies~\cite{petrakova2022mean,hota2023learning}.
In a congested power grid, an \gls{EV} with individual state of charge experiences different charging rates depending on the proportion of \gls{EV}'s charging at the same time~\cite{parise2014mean,bauso2017dynamic,tajeddini2018mean}.
Moreover, recent works have employed mean field games for the design of fair, non-monetary mechanisms in shared resource allocation settings, e.g., for traffic congestion management~\cite{elokda2022carma}.
In this case, the individual state is a private, stochastic urgency for the shared resource, as well as a scarce budget of tokens used to bid for the resource, whose dynamics depend on the macroscopic bidding behavior of others~\cite{balseiro2015repeated,elokda2023self}.

In all of these applications, the availability of practical computation tools to predict the \emph{\gls{MFE}} is instrumental to study the effect of complex control or policy design choices.
However, general formulations of mean field games face a computational bottleneck as they require solving an \emph{initial-terminal-value problem} consisting of a backward Hamilton-Jacobi-Bellman (HJB) equation characterizing the optimality of the agents’ policy trajectory; and a forward Fokker-Planck-Kolmogorov (FPK) equation ensuring that the mean field trajectory is consistent with the state dynamics under the optimal policy trajectory.
Existence (and uniqueness) of solutions to the initial-terminal-value problem is guaranteed under abstract assumptions that are difficult to verify in practice~\cite{huang2006large,lasry2007mean,gomes2010discrete,gomes2013continuous,almulla2017two,saldi2018markov}.
These assumptions typically manifest in linear dynamics with quadratic costs~\cite{parise2014mean,bauso2017dynamic,tajeddini2018mean,bensoussan2016linear}, and/or a weak dependency on the mean field (e.g., through the average state or action only~\cite{parise2014mean,tajeddini2018mean,aziz2016mean,moon2016linear}).

Due to these challenges, previous works have considered the less general solution concept of the \emph{stationary mean field equilibrium}~\cite{guo2023general,subramanian2019reinforcement,yardim2023policy}, also referred to as the \emph{\gls{SNE}}~\cite{adlakha2015equilibria,neumann2020stationary}, which is a single pair of optimal policy and the associated stationary mean field.
For the practically relevant class of \emph{discrete-time, finite-state-and-action, stationary mean-field games}, which we call \emph{\acrfull{DPGs}} for short, \gls{SNE} are computationally attractive as they are characterized as finite-dimensional fixed points rather than initial-terminal-value problems.
\gls{SNE} computation has been the subject of recent works, and thus far computational tools are based on standard Picard-type fixed point methods~\cite{guo2023general,yardim2023policy,cui2021approximately}.

A \emph{\gls{NE}} of a static \emph{population game} is also characterized by a fixed point;
however, there exists a rich set of \gls{NE} analysis and computation tools beyond conventional fixed point methods~\cite{sandholm2010population}.
In particular, \emph{evolutionary dynamics}, which include the well-known replicator, best response, Smith, and projection dynamics~\cite{sandholm2010population}, are useful for \gls{NE} computation,
since they leverage the game structure to ensure that stationary points coincide with \gls{NE}, and the stability of these dynamics is well-understood~\cite{sandholm2010population}.
Moreover, recent works have enriched the classical evolutionary dynamics with passivity-based and higher-order learning techniques that extend the classes of games in which convergence to a \gls{NE} is guaranteed~\cite{fox2013population,arcak2020dissipativity,gao2020passivity,park2021kl,toonsi2024higher}.

Population games are trivially mean field games, in which agents have only one state (with trivial dynamics) and finite actions.
The main contribution of this letter is to show that there is a larger and meaningful intersection: \acrfull{DPGs} \emph{are} population games.
We show a mathematical reduction of \gls{SNE} to standard \gls{NE} in population games:
For each \acrshort{DPG}, there exists a suitably defined static population game whose \gls{NE} coincide with \gls{SNE} of the \acrshort{DPG}.

The reduction of \gls{SNE} in \acrshort{DPGs} to \gls{NE} in population games unlocks a rich set of tools for \gls{SNE} analysis and computation.
We leverage it to guarantee the existence of a \gls{SNE}, develop an evolutionary dynamics-based \gls{SNE} computation algorithm, and provide simple conditions that guarantee stability and uniqueness of the \gls{SNE}.
Finally, the versatility of our approach is demonstrated in two complex application examples involving fair resource allocation to heterogeneous agents, as well as epidemic propagation and control.
Open source software for \gls{SNE} computation is provided at:\\
\url{https://gitlab.ethz.ch/elokdae/dynamic-population-games}.

The letter is organized as follows:
Section~\ref{sec:Preliminaries} defines \acrshort{DPGs};
Section~\ref{sec:ReductionPopGame} presents our main result: the reduction of \gls{SNE} to \gls{NE} in population games; Section~\ref{sec:ReductionConsequences} discusses important consequences of the reduction; and
Section~\ref{sec:applications} presents two application examples.
The letter is concluded in Section~\ref{sec:Conclusion}.

\subsection{Notation}
Let $a,d \in D \subseteq \Nat$  and let $c \in C\subseteq \Real^n$, then
for a function $f : D \times C \rightarrow \Real$, we distinguish discrete and continuous arguments through the notation $f[d](c)$.
Alternatively, we write $f : C \rightarrow \Real^{|D|}$ as the vector-valued function $f(c)$, with $f[d](c)$ denoting its $d$-{th} element.
Similarly, $p[a \mid d](c)$ denotes the conditional probability of $a$ given $d$ and $c$. 
For $m \in \Real_+$ denoting a measure of total mass, we denote by $\delta \in m \: \Delta(D):=\left\{\left. \sigma \in \Real_+^{|D|} \right\rvert \sum_{d \in D} \sigma[d] = m \right\}$ a distribution over the elements of $D$, with $\delta[d] \in [0, m]$ denoting the mass of element $d$ (if $m = 1$ then $\delta$ is a probability distribution).
Finally, when considering heterogeneous agent types, we denote by $x_\tau$ a quantity associated to type $\tau$.
\section{Dynamic Population Games (DPGs)}
\label{sec:Preliminaries}

In this section, we introduce the definition of a \emph{\gls{DPG}}, along with its solution concept, the \emph{\acrfull{SNE}}.
This class of games has been previously referred to as discrete-time, finite state-and-action, stationary mean-field games~\cite{guo2023general,subramanian2019reinforcement,yardim2023policy}.

\subsection{\acrfull{DPG} Setting}

A \gls{DPG}, denoted by $\G$, consists of a population \mbox{$\N = \{1,\dots,N\}$} of anonymous agents, with $N$ large such that the agents approximately form a \emph{continuum of mass}.
Each agent has one of a finite number of private static \emph{types} $\tau \in \T = \{1,\dots,N_\tau\}$, and the distribution of types is the exogenous parameter
$g \in \Delta(\T)$.
In addition to types, each agent has a private \emph{state} $x \in \X = \{1,\dots,N_x\}$,
which follows a dynamic process to be described hereafter. 
The \emph{state distribution} of type $\tau$ agents is given by $d_\tau \in \D_\tau = g_\tau \: \Delta(\X)$,
and the joint \emph{type-state distribution} is the concatenation of state distributions \mbox{$d=(d_1,\dots,d_{N_\tau}) \in \D = \prod_{\tau \in \T} \D_\tau$}.

At discrete time steps, agents participate in a game in which each agent of type-state $[\tau,x]$ plays one of a finite number of actions $a \in \A_{\tau}[x]$.
Type $\tau$ agents pick their action according to the stationary homogeneous \emph{policy} \mbox{$\pi_\tau : \X \rightarrow \Delta(\A_\tau[x])$},
which maps the state to a probability distribution over the actions\footnote{This class of policies has been previously referred to as \emph{oblivious policies}~\cite{adlakha2015equilibria}, since the agents choose their actions based on their individual state only and not the state distribution of others.}.
The probability to play action $a$ when in state $x$ is denoted by \mbox{$\pi_\tau[a \mid x] \in [0,1]$}. The set of policies of type $\tau$ is denoted by $\Pi_\tau$, and the concatenation of the policies is denoted by \mbox{$\pi=(\pi_1,\dots,\pi_{N_\tau}) \in \Pi = \prod_{\tau \in \T} \Pi_\tau$}.

The pair $(d,\pi) \in \D \times \Pi$ is referred to as the \emph{mean field} (in the tradition of mean field games) or the \emph{social state} (in the tradition of population games; a formal equivalency is shown in Section~\ref{sec:ReductionPopGame}).
It gives a macroscopic description of the distribution of the agents' types and states, as well as how they behave\footnote{Our definition of mean field $(d,\pi)$ differs slightly from the definition as the joint state-action distribution~\cite{jovanovic1988anonymous,adlakha2015equilibria,guo2023general}; however, the latter can be readily derived from $(d,\pi)$.}.

The discrete time state dynamics of type $\tau$ agents are governed by the \emph{state transition function}
\begin{align}
    \label{eq:Transitions}
    p_\tau[x^+ \mid x,a](d,\pi),
\end{align}
which is the probability of the state to transition from $x$ to $x^+$ after playing $a$, as a function of the mean field $(d,\pi)$.
The agent receives an \emph{immediate payoff} for playing $a$ when in state $x$, which is also a function of $(d,\pi)$ given by
\begin{align}
    \label{eq:Rewards}
    r_\tau[x,a](d,\pi).
\end{align}
A \gls{DPG} $\G$ is thus specified by the pair of state transition and immediate payoff functions~\eqref{eq:Transitions}--\eqref{eq:Rewards}.

\subsection{Solution Concept: \acrfull{SNE}}
\label{sec:StationaryNashEquilibrium}

Before introducing the solution concept in \gls{DPG}s, we must introduce definitions pertaining to the strategic problem faced by the population agents. 
Each agent faces a \emph{discounted \gls{MDP}}, i.e., aims to maximize its discounted infinite horizon payoff, with the type-dependant discount factor $\alpha_\tau \in [0,1)$.
Each agent's \gls{MDP} is coupled to other agents through the mean field $(d,\pi)$.
Namely, when a type $\tau$ agent follows policy $\pi_\tau$, its \emph{expected immediate payoff} and \emph{state transition matrix} are respectively given by,
\begin{align}
    R_\tau[x](d,\pi) &= \sum_{a \in \A_\tau[x]} \pi_\tau[a \mid x] \: r_\tau[x,a](d,\pi), \\
    P_\tau[x^+ \mid x](d,\pi) &= \sum_{a \in \A_\tau[x]} \pi_\tau[a \mid x] \: p_\tau[x^+ \mid x,a](d,\pi), \label{eq:StochasticMatrix}
\end{align}
and its \emph{expected infinite horizon payoff}, also known as the \emph{value function}, is recursively defined as 
\begin{multline}
    \label{eq:V-function-full}
    V_\tau[x](d,\pi) = R_\tau[x] + \alpha_\tau \sum_{x^+ \in \X} P_\tau[x^+ \mid x] \: V_\tau[x^+].
\end{multline}
Equation~\eqref{eq:V-function-full} is the Bellman recursion characterizing the value function associated with policy $\pi_\tau$.
Moreover, the \emph{single stage deviation payoff}, also known as the \emph{state-action value function} or \emph{$Q$-function}, is given by
\begin{multline}
    \label{eq:SingleStageDeviation}
    Q_\tau[x,a](d,\pi) = r_\tau[x,a] + \alpha_\tau \sum_{x^+ \in \X} p_\tau[x^+ \mid x,a] \: V_\tau[x^+].
\end{multline}
Equation~\eqref{eq:SingleStageDeviation} is the expected infinite horizon payoff if the type $\tau$ agent in state $x$ deviates from $\pi_\tau$ for a single stage by playing $a$, then follows $\pi_\tau$ in the future.
Notice that we omitted the dependency on $(d,\pi)$ on the right-hand side of Equations~\eqref{eq:V-function-full}--\eqref{eq:SingleStageDeviation} for notational convenience.

A policy attains the maximum value function~\eqref{eq:V-function-full} if and only if it maximizes the state-action value function~\eqref{eq:SingleStageDeviation} for all states $x \in \X$~\cite[Proposition~1.3.4]{bertsekas2007dynamic}.
This motivates the definition of the \emph{state-dependent best response correspondence}
{\small
\begin{multline}
    \label{eq:BestResponse}
    B_\tau[x](d,\pi) \in \\ \left\{\sigma \in \Delta(\A_\tau[x]) \left \lvert \begin{array}{l}
         \forall \sigma' \in \Delta(\A_\tau[x]), \\[2mm]
         \sum\limits_{a \in \A_\tau[x]} \left(\sigma[a] - \sigma'[a] \right) Q_\tau[x,a](d,\pi) \geq 0
    \end{array} \right. \right\}, \tag{BR}
\end{multline}
}
which is the set of probability distributions over the actions maximizing the state-action value function for a given type-state $[\tau,x] \in \T \times \X$.
With this, we are ready to define the solution concept in \gls{DPG}s.

\begin{definition}[\acrlong{SNE}]
\label{def:StationaryEquilibrium}
A \emph{\acrfull{SNE}} is a mean field $(\sd,\epi) \in \D \times \Pi$ which satisfies, for all $[\tau, x] \in \T \times \X$,
\begin{align}
    \sd_\tau[x] &= \sum_{x^- \in \X} \sd_\tau[x^-] \: P_\tau[x \mid x^-](\sd,\epi), \label{eq:SE-1} \tag{\gls{SNE}.1} \\
    \epi_\tau[\cdot \mid x] &\in B_{\tau}[x](\sd,\epi). \label{eq:SE-2} \tag{\gls{SNE}.2}
\end{align}
\end{definition}

The \gls{SNE} is similar to a classical \acrfull{NE} in that agents have no incentive to unilaterally deviate from the optimal policy $\epi$~\eqref{eq:SE-2}.
However, in addition, it requires that the type-state distribution $\sd$ is \emph{stationary} under the stochastic process characterized by $P_\tau(\sd,\epi)$~\eqref{eq:SE-1}.
The stationarity of $\sd$ ensures that the \gls{MDP} faced by each agent is time invariant at the equilibrium.

\section{Main Result: Reduction to Population Game}
\label{sec:ReductionPopGame}

In this section, we present the main result of this letter, which is a reduction of the \gls{SNE} in \gls{DPG}s to a standard \gls{NE} in static population games.
In order to state our result, we must first briefly recall the definition of population games~\cite{sandholm2010population}.
A population game is defined by a \emph{payoff vector} \mbox{$F(s)= (F_1(s),\dots,F_{N_\rho}(s))$}, where $F_\rho(s)$ is the payoff vector of population $\rho \in \Pop = \{1,\dots,N_\rho\}$.
Each entry of $F_\rho(s)$, denoted by $F_\rho[a](s)$, gives the payoff of a population $\rho$ agent for playing one of its finite actions $a \in \A_\rho$.
The \emph{social state} $s = (s_1,\dots,s_{N_\rho}) \in \Soc$ gives the distribution of actions in the population.
Namely, $s_\rho \in g_\rho \: \Delta(\A_\rho)$ is the distribution of actions in population $\rho$, where $g_\rho \in \Real_+$ is the total mass of agents in that population.
In the language of \gls{DPG}s, the populations $\rho$ can be viewed as the types $\tau$, and there are no time varying states $x$ (or, equivalently, the set of states $\X$ is a singleton).
The \emph{best response correspondence} of population $\rho$ is given by
{\small
\begin{multline}
    \label{eq:BestResponse-Pop}
    B_\rho(s) \in \\ \left\{\sigma \in g_\rho \: \Delta(\A_\rho) \left \lvert \begin{array}{l}
         \forall \sigma' \in g_\rho \: \Delta(\A_\rho), \\[2mm]
         \sum\limits_{a \in \A_\rho} \left(\sigma[a] - \sigma'[a] \right) F_\rho[a](s) \geq 0
    \end{array} \right. \right\}, \tag{PG-BR}
\end{multline}
}
and the \emph{\acrfull{NE}} is defined as a social state $\esoc$ that satisfies $\esoc_\rho \in B_\rho(\esoc)$ for all $\rho \in \Pop$.

\begin{theorem}[Reduction of \gls{SNE}]
\label{th:Reduction}
For every \gls{DPG} $\G$, there exists a static population game $F^\G$ whose Nash Equilibria (NE) coincide with the Stationary Nash Equilibria (SNE) of $\G$, i.e., it holds that $(\sd,\epi)$ is a \gls{SNE} of $\G$ if and only if $\esoc = (\sd,\epi)$ is a \gls{NE} of $F^\G$.
\end{theorem}

\begin{proof}
The proof proceeds by construction.
We construct a population game $F^\G$ with $N_\rho = N_\tau \: (1 + N_x)$ populations.
The first $N_\tau$ populations correspond to the different types, and we denote these populations as \mbox{$\rho = \tau \in \T$}.
The remaining $N_\tau N_x$ populations correspond to the \emph{type-state pairs}, and we denote these populations as \mbox{$\rho = [\tau,x] \in \T \times \X$}.
The action space of populations $\rho = \tau$ \emph{is} the state space $\X$, i.e., $\A_\tau = \X$, whereas the action space of populations $\rho = [\tau,x]$ is respectively \mbox{$\A_{[\tau,x]} = \A_\tau[x]$}.
The mass of populations $\rho = \tau$ is respectively $g_\tau$, whereas the mass of populations $\rho = [\tau,x]$ is uniformly $g_{[\tau,x]} = 1$.
With this construction, it holds that $\Soc = \D \times \Pi$, i.e., each social state $s \in \Soc$ under $F^\G$ is formally equivalent to a mean field $(d,\pi) \in \D \times \Pi$ under $\G$.
We thus write $s = (d,\pi)$ and define the payoffs $F^\G(d,\pi)$ for both groups of populations as
\begin{align}
    F^\G_\tau[x](d,\pi) &= \sum_{x^- \in \X} d_\tau[x^-] \: P_\tau[x \mid x^-](d,\pi) - d_\tau[x], \label{eq:payoff-dynamics} \\
    F^\G_{[\tau,x]}[a](d,\pi) &= Q_\tau[x,a](d,\pi). \label{eq:payoff-Q}
\end{align}
The payoffs of populations $\rho = \tau$ in Equation~\eqref{eq:payoff-dynamics} are motivated by condition~\eqref{eq:SE-1}.
Intuitively, these populations can be interpreted as \emph{agents of nature playing the dynamics}, i.e., they seek to bring the state distributions $d_\tau$ to stationarity.
On the other hand, the payoffs of populations $\rho = [\tau,x]$ in Equation~\eqref{eq:payoff-Q} coincide with the state-action value functions, as these are the quantities maximized by the agents in the \gls{DPG} $\G$.

To prove the theorem, we now show that the \gls{NE} condition of $F^\G$, i.e., $\esoc_\rho \in B_\rho(\esoc)$ for all $\rho \in \Pop$, is formally equivalent to conditions~\eqref{eq:SE-1}--\eqref{eq:SE-2}, with $\esoc = (\sd, \epi)$.
By construction, for all populations $\rho = [\tau,x] \in \T \times \X$, with \mbox{$\esoc_{[\tau,x]} = \epi_\tau[\cdot \mid x]$}, condition $\epi_\tau[\cdot \mid x] \in B_{[\tau,x]}(\sd, \epi)$ is formally equivalent to condition~\eqref{eq:SE-2}.
Therefore, we must show that for the remaining populations $\rho = \tau \in \T$, with $\esoc_\tau = \sd_\tau$, $\sd_\tau \in B_\tau(\sd,\epi)$ if and only if $\sd_\tau[x]$ satisfies condition~\eqref{eq:SE-1} for all $x \in \X$.
One direction of the implication is trivial: if $\sd_\tau[x]$ satisfies condition~\eqref{eq:SE-1} for all \mbox{$[\tau,x] \in \T \times \X$}, then \mbox{$F^\G_\tau[x](\sd,\epi) = 0$}, and all state distributions satisfy $d_\tau \in B_\tau(\sd,\epi)$, including $\sd_\tau$.
For the other direction of the implication, suppose, for the sake of contradiction, that $\sd_\tau \in B_\tau(\sd,\epi)$ for all $\tau \in \T$, but there is some type-state $[\tilde \tau, \tilde x] \in \T \times \X$ for which \mbox{$F^\G_{\tilde \tau}[\tilde x](\sd,\epi) \neq 0$}.
Notice that \mbox{$\sum_{x \in \X} F^\G_\tau[x](d,\pi) = 0$} for all $\tau \in \T$ and $(d,\pi) \in \D \times \X$; therefore, if $F^\G_{\tilde \tau}[\tilde x](\sd,\epi) > 0$ (respectively, $<0$), there must be another state \mbox{$\hat x \in \X$ for which $F^\G_{\tilde \tau}[\hat x](\sd,\epi) < 0$} (respectively, $>0$).
Suppose without loss of generality that $F^\G_{\tilde \tau}[\tilde x](\sd,\epi) > 0$ and $F^\G_{\tilde \tau}[\hat x](\sd,\epi) < 0$.
Since \mbox{$\sd_{\tilde \tau} \in B_{\tilde \tau}(\sd,\epi)$}, it must be that $\sd_{\tilde \tau}[\hat x] = 0$.
But this contradicts $F^\G_{\tilde \tau}[\hat x](\sd,\epi) < 0$, in which $-\sd_{\tilde \tau}[\hat x]$ is the only negative term.
We thus showed that if $\sd_\tau \in B_\tau(\sd,\epi)$ for all $\tau \in \T$, it must be that $F^\G_\tau[x](\sd,\epi) = 0$ for all \mbox{$[\tau,x] \in \T \times \X$}, i.e., $\sd_\tau$ satisfies condition~\eqref{eq:SE-1} for all \mbox{$[\tau,x] \in \T \times \X$}.
\end{proof}

\section{Consequences of the Reduction}
\label{sec:ReductionConsequences}

In this section, we discuss important consequences of our main result, Theorem~\ref{th:Reduction}, which establishes that \gls{SNE} in \gls{DPG}s can be reduced to standard \gls{NE} in static population games.
In general, all existing results pertaining to the \gls{NE} of the \emph{equivalent population game} $F^\G$, constructed in the proof of Theorem~\ref{th:Reduction}, can be adopted to the \gls{SNE} of \gls{DPG} $\G$.
We specialize a few of these results in this section.

\subsection{Existence of \gls{SNE}}

A \gls{SNE} is guaranteed to exist under the following mild assumption. This result can be proven directly using a fixed point argument~\cite{jovanovic1988anonymous,adlakha2015equilibria}, but in this letter we present a simpler proof based on Theorem~\ref{th:Reduction}.

\begin{assumption}[Continuity]
\label{ass:cont}
The state transition function $p_\tau(d,\pi)$ and the immediate payoff function $r_\tau(d,\pi)$ are continuous in the mean field $(d,\pi)$.
\end{assumption}

\begin{proposition}[Existence of \gls{SNE}]
    Let Assumption~\ref{ass:cont} hold for \gls{DPG} $\G$.
    Then $\G$ is guaranteed to have at least one \gls{SNE} $(\sd,\epi)$.
\end{proposition}

\begin{proof}
Consider the equivalent population game $F^\G(d,\pi)$ whose payoffs are given in Equations~\eqref{eq:payoff-dynamics}--\eqref{eq:payoff-Q}.
Assumption~\ref{ass:cont} implies that $F^\G(d,\pi)$ is continuous in $(d,\pi)$: continuity of~\eqref{eq:payoff-dynamics} is immediate; and continuity of~\eqref{eq:payoff-Q} follows from noticing that the value function~\eqref{eq:V-function-full} can be equivalently written in vector from as \mbox{$V_\tau(d,\pi) = \left(I - \alpha_\tau P_\tau(d,\pi)\right)^{-1} R_\tau(d,\pi)$}, which is continuous in $(d,\pi)$ for all $\alpha_\tau \in [0,1)$.
It follows that at least one \gls{NE} $(\sd,\epi)$ is guaranteed to exist in $F^\G$~\cite[Theorem~2.1.1]{sandholm2010population}, which coincides with a \gls{SNE} of $\G$ by Theorem~\ref{th:Reduction}.
\end{proof}

\subsection{Evolutionary Dynamics for \gls{SNE} Computation}
\label{sec:EvolutionaryPolStateDynamics}

Theorem~\ref{th:Reduction} moreover inspires the following \emph{evolutionary state-policy dynamics} model for \gls{SNE} computation
\begin{align}
    \dot{d}_\tau[x] &= \sum_{x^- \in \X} d_\tau[x^-] \: P_\tau[x \mid x^-](d,\pi) - d_\tau[x], \tag{EV.1} \label{eq:EvolutionState} \\
    \dot{\pi}_\tau[\cdot \mid x] &= \eta_\tau \: H_\tau(Q_\tau[x,\cdot](\pi,d),\pi_\tau[\cdot \mid x]), \tag{EV.2} \label{eq:EvolutionPolicy}
\end{align}
where $H_\tau$ is one of the familiar \emph{mean dynamics} in population games, such as the replicator, best response, Smith, or projection dynamic, whose closed-form expressions are given in~\cite[Chapters~4--6]{sandholm2010population}.
The rate parameter $\eta_\tau$ controls the speed of policy dynamics~\eqref{eq:EvolutionPolicy} with respect to state dynamics~\eqref{eq:EvolutionState}.
The dynamical system~\eqref{eq:EvolutionState}--\eqref{eq:EvolutionPolicy} is derived by applying evolutionary dynamics to the payoffs of the equivalent population game $F^\G(d,\pi)$ given in Equations~\eqref{eq:payoff-dynamics}--\eqref{eq:payoff-Q}.
In particular, the \emph{projection dynamic} is applied to~\eqref{eq:payoff-dynamics} to yield~\eqref{eq:EvolutionState}.
This choice is interpretable since it yields Laplacian diffusion dynamics in~\eqref{eq:EvolutionState}, however other mean dynamics could be used for the purpose of \gls{SNE} computation.
On the other hand, \eqref{eq:EvolutionPolicy} corresponds to applying mean dynamic $H_\tau$ to~\eqref{eq:payoff-Q}.
For computation purposes, it is natural to choose a mean dynamic that satisfies the well-known condition of \emph{Nash stationarity}~\cite[Section~5.2]{sandholm2010population} for $H_\tau$, to guarantee that stationary points of $H_\tau$ coincide with \gls{NE} of the underlying population game.
Most existing mean dynamics satisfy Nash stationarity, including the best response, Smith, and projection dynamic (but not the replicator dynamic)~\cite[Table~5.1]{sandholm2010population}.
The following proposition is then immediate from Theorem~\ref{th:Reduction}.

\begin{proposition}
    Suppose that mean dynamics $H_\tau$ satisfy Nash stationarity for all $\tau \in \T$.
    Then stationary points of the dynamical system~\eqref{eq:EvolutionState}--\eqref{eq:EvolutionPolicy} coincide with the \gls{SNE} of the \gls{DPG} $\G$.
\end{proposition}

A \gls{SNE} can thus be computed by forward simulating dynamical system~\eqref{eq:EvolutionState}--\eqref{eq:EvolutionPolicy}.
If the system reaches a stationary point, a \gls{SNE} is found.
Notice that
discretizing \eqref{eq:EvolutionState}--\eqref{eq:EvolutionPolicy} leads to many non-standard fixed point iterations of the \gls{SNE} map~\eqref{eq:SE-1}--\eqref{eq:SE-2} beyond the conventionally used Picard iteration~\cite{guo2023general,yardim2023policy,cui2021approximately}.
For example, the best response dynamic leads to a \emph{Krasnoselskij iteration} of~\eqref{eq:SE-1}--\eqref{eq:SE-2}, while other mean dynamics do not iterate the best response map~\eqref{eq:SE-2} at all.
The mean dynamic $H_\tau$, as well as the rate parameter $\eta_\tau$, are thus tuning parameters that can assist with the convergence to a \gls{SNE}.
For example, the Smith dynamic is better behaved than the best response dynamic in many examples of games~\cite{sandholm2010population}.
Moreover, (EV.1)--(EV.2) enable adopting state of the art, higher-order learning techniques that have recently extended the classes of games for which these dynamics converge~\cite{fox2013population,arcak2020dissipativity,gao2020passivity,park2021kl,toonsi2024higher}.

\subsection{Stability and Uniqueness of \gls{SNE}}

As a final consequence of Theorem~\ref{th:Reduction}, we introduce the following criterion for the stability and uniqueness of a \gls{SNE}.

\begin{definition}[Stable \gls{DPG}]
\label{def:stable-DPG}
\gls{DPG} $\G$ is a \emph{stable \gls{DPG}} if, for all $s=(d,\pi)$, $s'=(d',\pi') \in \D \times \Pi$,
{\small
\begin{multline}
    \label{eq:stable-DPG}
    \sum_{\tau \in \T} \left[\vphantom{\sum_{x \in \X}} \left(d_\tau - d'_\tau\right)^\top \left((P_\tau(s)^\top - I) \: d_\tau - (P_\tau(s')^\top - I) \: d'_\tau\right) + \right.\\
    \left. \sum_{x \in \X} \left(\pi_\tau[\cdot \mid x] - \pi'_\tau[\cdot \mid x]\right)^\top \left(Q_\tau[x,\cdot](s) - Q_\tau[x,\cdot](s')\right) \right] \leq 0.
\end{multline}
}
Moreover, if inequality~\eqref{eq:stable-DPG} holds strictly whenever $s' \neq s$, then $\G$ is a \emph{strictly stable \gls{DPG}}.
\end{definition}

Definition~\ref{def:stable-DPG} is an immediate application of the definition of \emph{stable population games}~\cite[Section~3.3]{sandholm2010population} to the equivalent population game $F^\G$.
Thus, the following proposition follows immediately from Theorem~\ref{th:Reduction} and~\cite[Theorem~3.3.10,~7.2.1,~7.2.4,~7.2.7,~7.2.9]{sandholm2010population}.

\begin{proposition}
    Let Assumption~\ref{ass:cont} hold for \gls{DPG} $\G$.
    If $\G$ is a \emph{stable \gls{DPG}}, then its \gls{SNE} are guaranteed to form a convex set.
    Moreover, the set of \gls{SNE} is globally asymptotically stable under the best response and Smith dynamics.\\
    If $\G$ is a \emph{strictly stable \gls{DPG}}, then it is guaranteed to have a unique \gls{SNE}, which is globally asymptotically stable under the replicator, best response, Smith and projection dynamics.
\end{proposition}
\section{Applications}
\label{sec:applications}

We present two application examples that demonstrate the versatility of our DPG formulation in tractably modelling complex real-world problems.
The modelling details of these examples are included in recent application-specific works that have utilized the general results of this letter~\cite{hota2023learning,elokda2022carma,elokda2023self,elokda2021adynamic,elokda2023dynamic,maitra2023sis}.

\subsection{Fair Resource Allocation}

Recent works have shown that introducing dynamics in a resource allocation problem enables the design of fair, non-monetary incentive schemes~\cite{elokda2022carma,elokda2023self,elokda2023dynamic}.
This is particularly important in contexts involving scarce, shared resources,
e.g., for traffic congestion management.
Instead of money, a non-tradable token, called \emph{karma}, is issued to the agents, who may use their karma to bid for shared resources.
Karma then flows from agents that acquire resources to those who yield resources, in accordance to exchange rules that are up to the policy maker to design.
The agents must be strategic and future sighted in their bidding, in order to spare karma for times of highest \emph{urgency}.
Thus, karma naturally induces a \gls{DPG} in which an agent's state is its current urgency and karma, and the immediate payoff for a given karma bid, as well as the karma dynamics, depend on the mean field through a potentially complex set of karma exchange rules.

Fortunately, the \gls{DPG} formulation has enabled devising tractable computational models of such a \emph{self-contained karma economy}.
In~\cite{elokda2023self}, computational insight on the \gls{SNE} of the karma \gls{DPG} is provided for different karma exchange rules, including first vs. second price auctions, peer-to-peer payments vs. payments based on societal redistribution, and even exchanges involving karma taxation.
Moreover, the fairness of the karma-based scheme with respect to heterogeneity in the agents' future awareness and urgency processes is investigated by incorporating multiple agent types.
In~\cite{elokda2022carma}, the self-contained karma economy is combined with a highway congestion model, where, in addition to bidding karma to access a congestion-free priority lane, agents must also choose when to travel in order to avoid the peak morning congestion, leading to even more complex and high-dimensional \gls{SNE}.
The model of~\cite{elokda2022carma} is further extended in~\cite{elokda2023dynamic}, where additionally an elastic demand of commuters is introduced as a new agent type.

\begin{figure}[!bt]
    \centering
    \includegraphics[height=.134\textwidth]{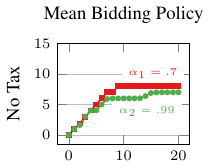}
    \hfill
    \includegraphics[height=.134\textwidth]{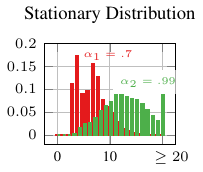}
    \hfill
    \includegraphics[height=.134\textwidth]{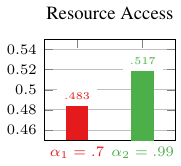}
    
    \includegraphics[height=.130\textwidth]{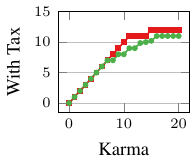}
    \hfill\hspace{0.5mm}
    \includegraphics[height=.130\textwidth]{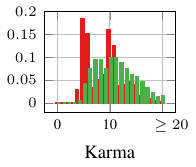}
    \hfill
    \includegraphics[height=.127\textwidth]{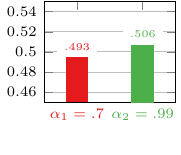}

    \caption{\gls{SNE} in karma \gls{DPG} with heterogeneous discount factors.}
    \label{fig:KarmaNash}

    

\end{figure}

In all of the above works, the SNE is computed using the evolutionary dynamics~\eqref{eq:EvolutionState}--\eqref{eq:EvolutionState}, with $H_\tau$ chosen as the \emph{perturbed best response dynamic}~\cite{sandholm2010population}.
Figure~\ref{fig:KarmaNash} shows two representative examples of karma \gls{SNE} adopted from~\cite{elokda2023self}, which correspond to a case where there are two agent types with different discount factors $\alpha_1=0.7$, $\alpha_2=0.99$.
The top and bottom row of Figure~\ref{fig:KarmaNash} allow to assess the effect of a small progressive karma tax (collected from each agent and redistributed) on the \gls{SNE}.
The left column of Figure~\ref{fig:KarmaNash} shows the mean bid of the \gls{SNE} policy $\epi$ ($y$-axis) per type and karma state ($x$-axis); the center column shows the stationary karma distribution $\sd$ per type; and the right column shows the stationary proportion of times each type gains access to the shared resource.
With no tax, the future sighted agents bid slightly less (Figure~\ref{fig:KarmaNash} top-left) and end up holding more karma on average (Figure~\ref{fig:KarmaNash} top-center) than the short sighted agents; as a consequence, the future sighted agents are better off in terms of average access to the resource (Figure~\ref{fig:KarmaNash} top-right).
The karma tax equalizes the field by redistributing karma to the short sighted agents.
To the extent of our knowledge, such complex and high-dimensional \gls{SNE} have not been computationally tractable before.

\subsection{Epidemic Modelling and Control}

\begin{figure}[!bt]
    \centering
    \includegraphics[height = .15\textwidth]{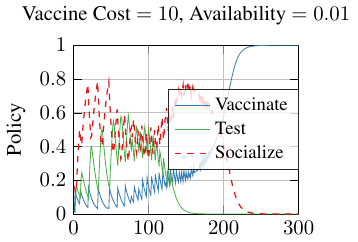}
    \hfil
    \includegraphics[height = .15\textwidth]{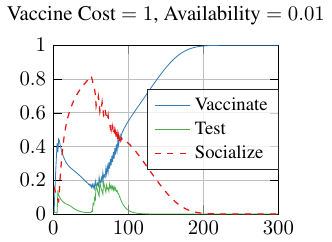}

    \hspace{1.5mm}
    \includegraphics[height = .154\textwidth]{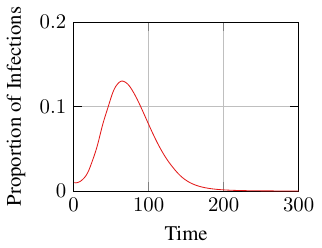}
    \hfill
    \includegraphics[height = .149\textwidth]{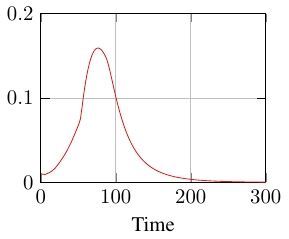}
    \hspace{5mm}
    
    \caption{Effect of vaccine cost and availability on epidemic spread.}
    \label{fig:Epidemic}
    
    
    
\end{figure}

Recent works have furthermore emphasized the importance of epidemic modelling and control as a mean field game, in which agents properly account for the future risk of infection rather than behave myopically~\cite{petrakova2022mean,hota2023learning,elokda2021adynamic,maitra2023sis}.
Among these works, \cite{hota2023learning,elokda2021adynamic,maitra2023sis} utilize the \gls{DPG} formulation to investigate the effect of complex policy interventions under minimal assumptions.
The evolutionary dynamics model~\eqref{eq:EvolutionState}--\eqref{eq:EvolutionState} is used not only for \gls{SNE} computation, but also as a model of the coupled evolution of infection states and population behaviors, which is needed to predict important off-equilibrium quantities such as the peak infections.
This computational approach, enabled by the results presented in this letter, allowed to study the effect of targeted lockdown restrictions in~\cite{elokda2021adynamic} and the interplay between cost and availability of vaccines and testing kits in~\cite{hota2023learning}. 
Figure~\ref{fig:Epidemic} shows a representative example adopted from~\cite{hota2023learning}, which demonstrates that decreasing the cost of vaccines without increasing its availability could lead to a surprising increase in peak infections (Figure~\ref{fig:Epidemic} right vs. left).
The top of Figure~\ref{fig:Epidemic} shows the policy evolution (i.e., agents' behavior): decreasing the cost of vaccines causes more asymptomatic agents to attempt to vaccinate rather than test for infection.
When they fail to do so due to limited vaccine availability, higher peak infections are observed (Figure~\ref{fig:Epidemic} bottom-right).

The availability of a tractable formulation and efficient numerical tools for this class of games allowed to perform extensive comparative studies of multiple interventions despite the complexity of this socio-technical system.


\section{Conclusion}
\label{sec:Conclusion}

We show that Stationary Nash Equilibria (SNE) in Dynamic Population Games (DPGs) can be reduced to standard Nash Equilibria (NE) in static population games.
This result unlocks a rich set of existing tools for SNE analysis and computation, which makes the DPG formulation favorable in comparison to other mean-field game formulations.
We demonstrate the versatility of the result in two complex application examples.
Important future topics include investigating the relationship between condition~(9) and previous conditions that guarantee convergence to a \gls{SNE}, as well as the decentralized computation of the \gls{SNE} by the agents.

\bibliographystyle{IEEEtran}
\bibliography{root}









\end{document}